\documentclass[11pt,reqno]{amsart}
\usepackage{amssymb,amscd,amsbsy}
\usepackage{amsthm}
\newcommand{\lb}{\linebreak}

\renewcommand{\a}{\alpha}
\renewcommand{\b}{\beta}
\newcommand{\g}{\gamma}

\newcommand{\e}{\varepsilon}

\newcommand{\z}{\zeta}

\renewcommand{\l}{\lambda}

\newcommand{\s}{\sigma}
\renewcommand{\t}{\tau}

\newcommand{\f}{\varphi}
\renewcommand{\o}{\omega}

\newcommand{\D}{\Delta}
\renewcommand{\L}{\Lambda}

\newcommand{\E}{{\mathcal E}}

\newcommand{\F}{{\mathcal F}}

\newcommand{\h}{{\mathcal H}}

\newcommand{\cS}{{\mathcal S}}
\newcommand{\X}{{\mathcal X}}
\newcommand{\Y}{{\mathcal Y}}

\newcommand{\T}{{\Bbb T}}
\newcommand{\pp}{{\Bbb P}}

\newcommand{\R}{{\Bbb R}}
\newcommand{\Z}{{\Bbb Z}}

\newcommand{\0}{{\boldsymbol{0}}}

\newcommand{\m}{{\boldsymbol m}}
\newcommand{\bS}{{\boldsymbol S}}

\newcommand{\cs}{\mbox{\rm const}\cdot}
\newcommand{\rf}[1]{(\ref{#1})}

\newcommand{\df}{\stackrel{\mathrm{def}}{=}}

\newcommand{\supp}{\operatorname{supp}}

\newcommand{\const}{\operatorname{const}}

\newcommand{\eeq}{\end{equation}}
\newcommand{\beq}{\begin{equation}}
\newcommand{\bay}{\begin{eqnarray}}
\newcommand{\ba}{\begin{align*}}
\newcommand{\ea}{\end{align*}}
\newcommand{\ey}{\end{eqnarray}}
\newcommand{\bey}{\begin{eqnarray*}}
\newcommand{\eey}{\end{eqnarray*}}

\newcommand{\imp}{\Rightarrow}
\newcommand{\be}{\infty}

\newcommand{\bl}{\blacksquare}

\newcommand{\Pf}{{\bf Proof. }}

\newcommand{\ov}{\overline}

\newtheorem{thm}{\hspace{\parindent}Theorem}[section]

\newtheorem{lem}[thm]{\hspace{\parindent}Lemma}

\theoremstyle{remark}

\newtheorem*{rem*}{Remark}

\newcommand\fM{\frak M}
\newcommand\cZ{\mathcal{Z}}
\newcommand\dg{\frak D}
\newcommand\pt{\hat\otimes_{\rm i}}

\begin{document}

\newcommand{\vse}{\vspace{.2in}}
\numberwithin{equation}{section}

\title{\bf Multiple operator integrals \\ and higher operator derivatives}
\author{V.V. Peller}
\thanks{The author is partially supported by NSF grant DMS 0200712}
\maketitle

\newcommand{\mt}{{\mathcal T}}

\begin{abstract}
In this paper we consider the problem of the existence of higher derivatives
of the function $t\mapsto\f(A+tK)$, where $\f$ is a function on the real line,
$A$ is a self-adjoint operator, and $K$ is a bounded self-adjoint operator.
We improve earlier results by Sten'kin. In order to do this, we give a new approach 
to multiple operator integrals. This approach improves the earlier approach given by Sten'kin.
We also consider a similar problem for unitary operators.
\end{abstract}

\setcounter{section}{0}
\section{\bf Introduction}
\setcounter{equation}{0}

\

If $A$ is a bounded self-adjoint operator on Hilbert space, the spectral theorem allows one
for Borel function $\f$ on the real line $\R$ to define the function $\f(A)$ of $A$.
We are going to study in this paper smoothness properties of the map $A\mapsto\f(A)$.
It is easy to see that if this map is differentiable (in the sense of Gateaux), then
$\f$ is continuously differentiable.

If $K$ is another bounded self-adjoint operator, consider the function $t\mapsto\f(A+tK)$,
$t\in\R$. In \cite{DK} it was shown that if $\f\in C^2(\R)$ (i.e., $\f$ is twice continuously
differentiable), then the map $t\mapsto\f(A+tK)$ is norm differentiable and
\bay
\label{DKF}
\frac{d}{ds}\big(\f(A+sK)\big)\Big|_{s=0}=\iint\frac{\f(\l)-\f(\mu)}{\l-\mu}\,dE_A(\l)K\,dE_A(\mu),
\ey
where $E_A$ is the spectral measure of $A$. The expression on the right-hand side of \rf{DKF}
is a {\it double operator integral}. Later Birman and Solomyak developed their beautiful theory of
double operator integrals in \cite{BS1}, \cite{BS2}, and \cite{BS3} (see also \cite{BS4});
we discuss briefly this theory in \S 3.

Throughout this paper {\it if we integrate a function on $\R^d$ (or $\T^d$) and the domain of
integration is not specified, it is assumed that the domain of integration is $\R^d$ (or $\T^d$).}

Birman and Solomyak relaxed in \cite{BS3} the assumptions on $\f$ under which \rf{DKF} holds.
They also considered the case of an unbounded self-adjoint operator $A$. However, it turned out that
the condition $\f\in C^1(\R)$ is not sufficient for the differentiability  of the function
$t\mapsto\f(A+tK)$ even in the case of bounded $A$. This can be deduced from an explicit example
constructed by Farforovskaya \cite{F2} (in fact, this can also be deduced from an example given
in \cite{F1}). 

In \cite{Pe1} a necessary condition on $\f$ for the differentiability of the function 
\lb$t\mapsto\f(A+tK)$ for all $A$ and $K$ was found. That necessary condition was deduced from the
nuclearity criterion for Hankel operators (see the monograph \cite{Pe4}) and it implies that
the condition $\f\in C^1(\R)$ is not sufficient. We also refer the reader to \cite{Pe2} where a
necessary condition is given in the case of an unbounded self-adjoint operator $A$. 

Sharp sufficient conditions on $\f$ for the differentiability of the function 
$t\mapsto\f(A+tK)$ were obtained in \cite{Pe1} in the case of bounded self-adjoint operators 
and in \cite{Pe2} in the case of an unbounded self-adjoint operator $A$. In particular, it follows
from the results of \cite{Pe2} that if $\f$ belongs to the homogeneous Besov space 
$B_{\be1}^1(\R),$\footnote{see \S 2 for information on Besov spaces}
$A$ is a self-adjoint operator and $K$ is a bounded self-adjoint operator, then the function
$t\mapsto\f(A+tK)$ is differentiable and \rf{DKF} holds (see \S 5 of this paper for details).
In the case of bounded self-adjoint operators formula \rf{DKF} holds if $\f$ belongs to 
$B_{\be1}^1(\R)$ locally (see \cite{Pe1}).

A similar problem for unitary operators was considered in \cite{BS3} and
later in \cite{Pe1}. Let $\f$ be a function on the unit circle $\T$.
For a unitary operator $U$ and a bounded self-adjoint operator $A$, consider the function
$t\mapsto\f(e^{{\rm i}tA}U)$. It was shown in \cite{Pe1} that if $\f$ belongs to the Besov space 
$B_{\be1}^1$, then the function $t\mapsto\f(e^{{\rm i}tA}U)$ is differentiable and
\bay
\label{du}
\frac{d}{ds}\Big(\f(e^{{\rm i}sA}U)\Big)\Big|_{s=o}
={\rm i}\left(\iint\frac{\f(\l)-\f(\mu)}{\l-\mu}\,dE_U(\l)A\,dE_U(\mu)\right)U
\ey
(earlier this formula was obtained in \cite{BS3} under more restrictive assumptions on $\f$). We refer
the reader to \cite{Pe1} and \cite{Pe2} for necessary conditions. We also mention here the paper 
\cite {ABF}, which slightly improves the sufficient condition
$\f\in B_{\be1}^1$.

The problem of the existence of higher derivatives of the function $t\mapsto\f(A+tK)$
was studied by Sten'kin in \cite{S}. He showed that under certain conditions on $\f$ the function
$t\mapsto\f(A+tK)$ has $m$ derivatives and
\bay
\label{SF}
\frac{d^m}{ds^m}\big(\f(A_s)\big)\Big|_{s=0}\!\!=
m!\underbrace{\int\!\!\cdots\!\!\int}_{m}(\dg^{m}\f)(\l_1,\cdots,\l_{m+1})
\,dE_A(\l_1)K\cdots K\,dE_A(\l_{m+1}),
\ey
where the {\it divided differences $\dg^k\f$ of order $k$} are defined inductively by
$$
\dg^0\f\df\f,
$$
$$
(\dg^{k}\f)(\l_1,\cdots,\l_{k+1})\df
\frac{(\dg^{k-1}\f)(\l_1,\cdots,\l_k)-(\dg^{k-1}\f)(\l_2,\cdots,\l_{k+1})}{\l_{1}-\l_{k+1}},\quad k\ge1
$$
(the definition does not depend on the order of the variables). We are also going to use
the notation
$$
\dg\f=\dg^1\f.
$$

The Birman--Solomyak theory of double
operator integrals does not generalize to the case of multiple operator integrals. In \cite{S} 
to define the multiple operator integrals 
$$
\underbrace{\int\cdots\int}_k\psi(\l_1,\cdots,\l_k)
\,dE_1(\l_1)T_1\,dE_2(\l_2)T_2\cdots T_{k-1}\,dE_k(\l_k),
$$
Sten'kin considered iterated integration and he defined such integrals
for a certain class of functions $\psi$. However, that approach in the case $k=2$ leads to a
considerably smaller class of functions $\psi$ than the Birman--Solomyak approach.
In particular the function $\psi$ identically equal to $1$, is not integrable in the sense of
the approach developed in \cite{S}, while it is very natural to assume that
$$
\underbrace{\int\cdots\int}_k
\,dE_1(\l_1)T_1\,dE_2(\l_2)T_2\cdots T_{k-1}\,dE_k(\l_k)=T_1T_2\cdots T_{k-1}.
$$  

In \S 3 of this paper we use a different approach to the definition of multiple operator integrals.
The approach is based on integral projective tensor products. In the case $k=2$ our approach
produces the class of integrable functions that coincides with the class of so-called
Schur multipliers, which is the maximal possible class of integrable functions
in the case $k=2$ (see \S 4).

Our approach allows us to improve Sten'kin's results on the existence of higher order derivatives of
the function $t\mapsto\f(A+tK)$. We prove in \S 5 that formula \rf{SF} holds for functions
$\f$ in the intersection $B_{\be1}^m(\R)\bigcap B_{\be1}^1(\R)$ of homogeneous Besov spaces.

Note that the Besov spaces $B_{\be1}^1$ and $B_{\be1}^1(\R)$ appear in a natural way when
studying the applicability of the Lifshits--Krein trace formula for trace class perturbations
(see \cite{Pe1} and \cite{Pe2}), while the Besov spaces $B_{\be1}^2$ and $B_{\be1}^2(\R)$ arise
when studying the applicability of the Koplienko--Neidhardt trace formulae for Hilbert--Schmidt 
perturbations (see \cite{Pe5}).

It is also interesting to note that the Besov class $B_{\be1}^2(\R)$ appears in a natural way
in perturbation theory in \cite{Pe3}, where the following problem is studied: in which case
$$
\f(T_f)-T_{\f\circ f}\in \bS_1?
$$ 
($T_g$ is a Toeplitz operator with symbol $g$.)

In \S 4 we obtain similar results in the case of unitary operators and generalize formula
\rf{du} to the case of higher derivatives.

We collect in \S 3 basic information on Besov spaces.

\

\section{\bf Besov spaces}
\setcounter{equation}{0}

\

 Let $0<p,\,q\le\be$ and $s\in\R$. The Besov class $B^s_{pq}$ of functions (or
distributions) on $\T$ can be defined in the following way. Let $w$ be a $C^\be$ function on $\R$ such
that
\bay
\label{v}
w\ge0,\quad\supp w\subset\left[\frac12,2\right],\quad\mbox{and}\quad
\sum_{n=-\be}^\be w(2^{n}x)=1\quad\mbox{for}\quad x>0.
\ey
Consider the trigonometric polynomials $W_n$, and $W_n^\#$ defined by
$$
W_n(z)=\sum_{k\in\Z}w\left(\frac{k}{2^n}\right)z^k,\quad n\ge1,\quad W_0(z)=\bar z+1+z,\quad
\mbox{and}\quad W_n^\#(z)=\ov{W_n(z)},\quad n\ge0.
$$
Then for each distribution $\f$ on $\T$,
$$
\f=\sum_{n\ge0}\f*W_n+\sum_{n\ge1}\f*W^\#_n.
$$
The Besov class $B^s_{pq}$ consists of functions (in the case $s>0$) or distributions $\f$ on $\T$
such that
$$
\big\{\|2^{ns}\f*W_n\|_{L^p}\big\}_{n\ge0}\in\ell^q\quad\mbox{and}
\quad\big\{\|2^{ns}\f*W^\#_n\|_{L^p}\big\}_{n\ge1}\in\ell^q
$$
Besov classes admit many other descriptions. In particular, for $s>0$, the space $B^s_{pq}$ admits the
following characterization. A function $\f$ belongs to $B^s_{pq}$, $s>0$, if and only if
$$
\int_\T\frac{\|\D^n_\t f\|_{L^p}^q}{|1-\t|^{1+sq}}d\m(\t)<\be\quad\mbox{for}\quad q<\be
$$
and
$$
\sup_{\t\ne1}\frac{\|\D^n_\t f\|_{L^p}}{|1-\t|^s}<\be\quad\mbox{for}\quad q=\be,
$$
where $\m$ is normalized Lebesgue measure on $\T$, $n$ is an integer greater than $q$ and $\D_\t$ is the
difference operator: 
$(\D_\t f)(\z)=f(\t\z)-f(\z)$, $\z\in\T$.

To define (homogeneous) Besov classes $B^s_{pq}(\R)$ on the real line, we consider the same function $w$
as in \rf{v} and define the functions $W_n$ and $W^\#_n$ on $\R$ by
$$
\F W_n(x)=w\left(\frac{x}{2^n}\right),\quad\F W^\#_n(x)=\F W_n(-x),\quad n\in\Z,
$$
where $\F$ is the {\it Fourier transform}. The Besov class $B^s_{pq}(\R)$ consists of
distributions $\f$ on $\R$ such that
$$
\{\|2^{ns}\f*W_n\|_{L^p}\}_{n\in\Z}\in\ell^q(\Z)\quad\mbox{and}
\quad\{\|2^{ns}\f*W^\#_n\|_{L^p}\}_{n\in\Z}\in\ell^q(\Z).
$$
According to this definition, the space $B^s_{pq}(\R)$ contains all polynomials. However, it is not
necessary to include all polynomials. 

In this paper we need only Besov spaces $B_{\be1}^d$, $d\in\Z_+$. In the case of functions
on the real line it is convenient to restrict the degree of polynomials in $B_{\be1}^d(\R)$ by $d$.
It is also convenient to consider the following seminorm on
$B_{\be1}^d(\R)$:
$$
\|\f\|_{B_{\be1}^d(\R)}=\sup_{x\in\R}|\f^{(d)}(x)|+\sum_{n\in\Z}2^{nd}\|\f*W_n\|_{L^\be}+
\sum_{n\in\Z}2^{nd}\|\f*W^\#_n\|_{L^\be}.
$$
The classes $B_{\be1}^d(\R)$ can be described as classes of function on $\R$ 
in the following way:
$$
\f\in B_{\be1}^d(\R)\quad\Longleftrightarrow\quad
\sup_{t\in\R}|\f^{(d)}(t)|+\int\limits_\R\frac{\|\D_t^{d+1}\f\|_{L^\be}}{|t|^{d+1}}\,dt<\be,
$$
where $\D_t$ is the difference operator defined by $(\D_t\f)(x)=\f(x+t)-\f(x)$.

We refer the reader to \cite{Pee} for more detailed information on Besov classes.

\

\section{\bf Multiple operator integrals}
\setcounter{equation}{0}

\

In this section we define multiple operator integrals using integral projective tensor products of
$L^\be$-spaces. However, we begin with a brief review of the theory of double operator integrals
that was developed by Birman and Solomyak in \cite{BS1}, \cite{BS2}, and \cite{BS3}. We state a
description of the Schur multipliers associated with two spectral measures in terms of integral
projective tensor products. This suggests the idea to define multiple operator integrals with the
help of integral projective tensor products.

\medskip

{\bf Double operator integrals.} Let $(\X,E)$ and $(\Y,F)$ be spaces with spectral measures $E$ and $F$
on a Hilbert space $\h$. Let us first define double operator integrals
\bay
\label{doi}
\int\limits_\X\int\limits_\Y\psi(\l,\mu)\,d E(\l)T\,dF(\mu),
\ey
for bounded measurable functions $\psi$ and operators $T$
of Hilbert Schmidt class $\bS_2$. Consider the spectral measure $\E$ whose values are orthogonal
projections on the Hilbert space $\bS_2$, which is defined by
$$
\E(\L\times\D)T=E(\L)TF(\D),\quad T\in\bS_2,
$$ 
$\L$ and $\D$ being measurable subsets of $\X$ and $\Y$. Then $\E$ extends to a spectral measure on 
$\X\times\Y$ and if $\psi$ is a bounded measurable function on $\X\times\Y$, by definition,
$$
\int\limits_\X\int\limits_\Y\psi(\l,\mu)\,d E(\l)T\,dF(\mu)=
\left(\,\,\int\limits_{\X\times\Y}\psi\,d\E\right)T.
$$
Clearly,
$$
\left\|\int\limits_\X\int\limits_\Y\psi(\l,\mu)\,dE(\l)T\,dF(\mu)\right\|_{\bS_2}
\le\|\psi\|_{L^\be}\|T\|_{\bS_2}.
$$
If 
$$
\int\limits_\X\int\limits_\Y\psi(\l,\mu)\,d E(\l)T\,dF(\mu)\in\bS_1
$$
for every $T\in\bS_1$, we say that $\psi$ is a {\it Schur multiplier (of $\bS_1$) associated with 
the spectral measure $E$ and $F$}. In
this case by duality the map
\bay
\label{tra}
T\mapsto\int\limits_\X\int\limits_\Y\psi(\l,\mu)\,d E(\l)T\,dF(\mu),\quad T\in \bS_2,
\ey
extends to a bounded linear transformer on the space of bounded linear operators on $\h$.
We denote by $\fM(E,F)$ the space os Schur multipliers of $\bS_1$ associated with
the spectral measures $E$ and $F$. The norm of $\psi$ in $\fM(E,F)$ is, by definition, the norm of the
transformer \rf{tra} on the space of bounded linear operators.

In \cite{BS3} it was shown that if $A$ is a self-adjoint operator (not necessarily bounded),
$K$ is a bounded self-adjoint operator and if
$\f$ is a continuously differentiable 
function on $\R$ such that the divided difference $\dg\f$ is a Schur multiplier
of $\bS_1$ with respect to the spectral measures of $A$ and $A+K$, then
\bay
\label{BSF}
\f(A+K)-\f(A)=\iint\frac{\f(\l)-\f(\mu)}{\l-\mu}\,dE_{A+K}(\l)K\,dE_A(\mu)
\ey
and
$$
\|\f(A+K)-\f(A)\|\le\const\|\f\|_{\fM(E_A,E_{A+K})}\|K\|,
$$
i.e., {\it $\f$ is an operator Lipschitz function}.

It is easy to see that if a function $\psi$ on $\X\times\Y$ belongs to the {\it projective tensor
product}
$L^\be(E)\hat\otimes L^\be(F)$ of $L^\be(E)$ and $L^\be(F)$ (i.e., $\psi$ admits a representation
$$
\psi(\l,\mu)=\sum_{n\ge0}f_n(\l)g_n(\mu),
$$
where $f_n\in L^\be(E)$, $g_n\in L^\be(F)$, and
$$
\sum_{n\ge0}\|f_n\|_{L^\be}\|g_n\|_{L^\be}<\be),
$$
then $\psi\in\fM(E,F)$.
For such functions $\psi$ we have
$$
\int\limits_\X\int\limits_\Y\psi(\l,\mu)\,d E(\l)T\,dF(\mu)=
\sum_{n\ge0}\left(\,\int\limits_\X f_n\,dE\right)T\left(\,\int\limits_\Y g_n\,dF\right).
$$ 

More generally, $\psi$ is a Schur multiplier of $\bS_1$ if $\psi$ 
belongs to the {\it integral projective tensor product} $L^\be(E)\hat\otimes_{\rm i}
L^\be(F)$ of $L^\be(E)$ and $L^\be(F)$, i.e., $\psi$ admits a representation
\bay
\label{ipt}
\psi(\l,\mu)=\int_Q f(\l,x)g(\mu,x)\,d\s(x),
\ey
where $(Q,\s)$ is a measure space, $f$ is a measurable function on $\X\times Q$,
$g$ is a measurable function on $\Y\times Q$, and
\bay
\label{ir}
\int_Q\|f(\cdot,x)\|_{L^\be(E)}\|g(\cdot,x)\|_{L^\be(F)}\,d\s(x)<\be.
\ey
If $\psi\in L^\be(E)\hat\otimes_{\rm i}L^\be(F)$, then
$$
\int\limits_\X\int\limits_\Y\psi(\l,\mu)\,d E(\l)T\,dF(\mu)=
\int\limits_Q\left(\,\int\limits_\X f(\l,x)\,dE(\l)\right)T
\left(\,\int\limits_\Y g(\mu,x)\,dF(\mu)\right)\,d\s(x).
$$

It turns out that all Schur multipliers can be obtained in this way. More precisely, the following
result holds (see \cite{Pe1}):

\medskip

{\bf Theorem on Schur multipliers.} {\em Let $\psi$ be a measurable function on 
$\X\times\Y$. The following are equivalent

{\rm (i)} $\fM(E,F)$;

{\rm (ii)} $\psi\in L^\be(E)\hat\otimes_{\rm i}L^\be(F)$;

{\rm (iii)} there exist measurable functions $f$ on $\X\times Q$ and $g$ on $\Y\times Q$ such that
{\em\rf{ipt}} holds and
\bay
\label{bs}
\left\|\int_Q|f(\cdot,x)|^2\,d\s(x)\right\|_{L^\be(E)}
\left\|\int_Q|g(\cdot,x)|^2\,d\s(x)\right\|_{L^\be(F)}<\be.
\ey
}

Note that the implication (iii)$\imp$(ii) was established in \cite{BS3}. Note also that
in the case of matrix Schur multipliers (this corresponds to discrete spectral measures
of multiplicity 1) the equivalence of (i) and (ii) was proved in \cite{Be}.

It is interesting to observe that if $f$ and $g$ satisfy \rf{ir}, then they also satisfy
\rf{bs}, but the converse is false. However, if $\psi$ admits a representation of the form \rf{ipt}
with $f$ and $g$ satisfying \rf{bs}, then it also admits a (possibly different) representation of the
form \rf{ipt} with $f$ and $g$ satisfying \rf{ir}.

Note that in a similar way we can define the {\it projective tensor product} $A\hat\otimes B$
and the {\it integral projective tensor product} $A\hat\otimes_{\rm i} B$ of 
arbitrary Banach functions spaces $A$ and $B$.

The equivalence of (i) and (ii) in the Theorem on Schur multipliers suggests an idea how to define
multiple operator integrals.

\medskip

{\bf Multiple operator integrals.} 
We can easily extend the definition of the projective tensor product and the integral
projective tensor product to three or more function spaces.

Consider first the case of triple operator integrals.

Let $(\X,E)$, $(\Y,F)$, and $(\cZ,G)$
be spaces with spectral measures $E$, $F$, and $G$ on a Hilbert space $\h$. Suppose that
$\psi$ belongs to the integral projective tensor product
$L^\be(E)\hat\otimes_{\rm i}L^\be(F)\hat\otimes_{\rm i}L^\be(G)$, i.e., $\psi$ admits a representation
\bay
\label{ttp}
\psi(\l,\mu,\nu)=\int_Q f(\l,x)g(\mu,x)h(\nu,x)\,d\s(x),
\ey
where $(Q,\s)$ is a measure space, $f$ is a measurable function on $\X\times Q$,
$g$ is a measurable function on $\Y\times Q$, $h$ is a measurable function on $\cZ\times Q$,
and
\bay
\label{ner}
\int_Q\|f(\cdot,x)\|_{L^\be(E)}\|g(\cdot,x)\|_{L^\be(F)}\|h(\cdot,x)\|_{L^\be(G)}\,d\s(x)<\be.
\ey
We define the norm $\|\psi\|_{L^\be\hat\otimes_{\rm i}L^\be\hat\otimes_{\rm i}L^\be}$ in the
space $L^\be(E)\hat\otimes_{\rm i}L^\be(F)\hat\otimes_{\rm i}L^\be(G)$ as the infimum
of the left-hand side of \rf{ner} over all representations 
\rf{ttp}. 

Suppose now that $T_1$ and $T_2$ be bounded linear operators on $\h$. For a function $\psi$ in
$L^\be(E)\hat\otimes_{\rm i}L^\be(F)\hat\otimes_{\rm i}L^\be(G)$ of the form \rf{ttp}, we put
\begin{align}
\label{opr}
&\int\limits_\X\int\limits_\Y\int\limits_\cZ\psi(\l,\mu,\nu)
\,d E(\l)T_1\,dF(\mu)T_2\,dG(\nu)\nonumber\\[.2cm]
\df&\int\limits_Q\left(\,\int\limits_\X f(\l,x)\,dE(\l)\right)T_1
\left(\,\int\limits_\Y g(\mu,x)\,dF(\mu)\right)T_2
\left(\,\int\limits_\cZ h(\nu,x)\,dG(\nu)\right)\,d\s(x).
\end{align}

The following lemma shows that the triple operator integral
$$
\int\limits_\X\int\limits_\Y\int\limits_\cZ\psi(\l,\mu,\nu)\,d E(\l)T_1\,dF(\mu)T_2\,dG(\nu)
$$
is well-defined.

\begin{lem}
\label{kor}
Suppose that $\psi\in L^\be(E)\hat\otimes_{\rm i}L^\be(F)\hat\otimes_{\rm i}L^\be(G)$. Then the 
right-hand side of {\em\rf{opr}} does not depend on the choice of a representation {\em\rf{ttp}}
and
\bay
\label{in}
\left\|\int\limits_\X\int\limits_\Y\int\limits_\cZ\psi(\l,\mu,\nu)
\,dE(\l)T_1\,dF(\mu)T_2\,dG(\nu)\right\|
\le\|\psi\|_{L^\be\hat\otimes_{\rm i}L^\be\hat\otimes_{\rm i}L^\be}\cdot\|T_1\|\cdot\|T_2\|.
\ey
\end{lem}

\Pf To show that the right-hand side of \rf{opr} does not depend on the choice of a representation
\rf{ttp}, it suffices to show that if the right-hand side of \rf{ttp} is the zero function,
then the right-hand side of \rf{opr} is the zero operator. Denote our Hilbert space by $\h$
and let $\z\in\h$. We have
$$
\int\limits_\cZ\left(\int\limits_Q f(\l,x)g(\mu,x)h(\nu,x)\,d\s(x)\right)\,dG(\nu)=0\quad
\mbox{for almost all $\l$ and $\mu$},
$$
and so for almost all $\l$ and $\mu$,
\begin{align*}
&\int\limits_Q f(\l,x)g(\mu,x)T_2\left(\,\int\limits_\cZ h(\nu,x)\,dG(\nu)\right)\z\,d\s(x)\\=
&T_2\int\limits_\cZ\left(\int\limits_Q f(\l,x)g(\mu,x)h(\nu,x)\,d\s(x)\right)\,dG(\nu)\z=\0.
\end{align*}
Putting
$$
\xi_x=T_2\left(\,\int\limits_\cZ h(\nu,x)\,dG(\nu)\right)\z,
$$
we obtain
$$
\int\limits_Q f(\l,x)g(\mu,x)\xi_x\,d\s(x)=\0\quad\mbox{for almost all}\quad \l\quad\mbox{and}\quad
\mu.
$$
We can realize the Hilbert space $\h$ as a space of vector functions so that integration
with respect to the spectral measure $F$ corresponds to multiplication. It follows that
$$
\int\limits_Q \!f(\l,x)T_1\!\left(\int\limits_\Y \!g(\mu,x)\,dF(y)\!\right)\!\xi_x\,d\s(x)\,dF(\mu)=
T_1\!\int\limits_\Y\!\int\limits_Q \!f(\l,x)g(\mu,x)\xi_x\,d\s(x)\,dF(\mu)=\0
$$
for almost all $\l$. Let now
$$
\eta_x=T_1\left(\int\limits_\Y \!g(\mu,x)\,dF(\mu)\right)\xi_x.
$$
We have
$$
\int\limits_Q f(\l,s)\eta_x\,d\s(x)=\0\quad\mbox{for almost all}\quad\l.
$$
Now we can realize $\h$ as a space of vector functions so that integration
with respect to the spectral measure $E$ corresponds to multiplication. It follows that
$$
\int\limits_Q \left(\int\limits_\X f(\l,x)\,dE(\l)\right)\eta_x\,d\s(x)=
\int\limits_\X\int\limits_Q f(\l,x)\eta_x\,d\s(x)\,dE(\l)=\0.
$$
This exactly means that the right-hand side of \rf{opr} is the zero operator.

Inequality \rf{in} follows immediately from \rf{opr}. $\bl$

In a similar way we can define multiple operator integrals
$$
\underbrace{\int\cdots\int}_k\psi(\l_1,\cdots,\l_k)
\,dE_1(\l_1)T_1\,dE_2(\l_2)T_2\cdots T_{k-1}\,dE_k(\l_k)
$$
for functions $\psi$ in the integral projective tensor product 
$\underbrace{L^\be(E_1)\hat\otimes_{\rm i}\cdots\hat\otimes_{\rm i}L^\be(E_k)}_k$
(the latter space is defined in the same way as in the case $k=2$).

\

\section{\bf The case of unitary operators}
\setcounter{equation}{0}

\

Let $U$ be a unitary operator and $A$ a bounded
self-adjoint on Hilbert space. For $t\in\R$, we put
$$
U_t=e^{{\rm i}tA}U.
$$
In this section we obtain sharp conditions on the existence of higher operator derivatives
of the function $t\mapsto \f(U_t)$.

Recall that it was proved in \cite{Pe1} that for a function $\f$ in the Besov space $B_{\be1}^1$
the divided difference $\dg\f$ belongs to the projective tensor product
$C(\T)\hat\otimes C(\T)$, and so for arbitrary unitary operators $U$ and $V$
the following formula holds:
\bay
\label{vp}
\f(V)-\f(U)=\iint\frac{\f(\l)-\f(\mu)}{\l-\mu}\,dE_V(\l)(V-U)\,dE_U(\mu).
\ey

First we state the main results of this section for second derivatives.

\begin{thm}
\label{tens}
If $\f\in B_{\be1}^2$, then 
$$
(\dg^2\f)\in C(\T)\hat\otimes C(\T)\hat\otimes C(\T).
$$
\end{thm}

\begin{thm}
\label{uni}
Let $\f$ be a function in the Besov class $B^2_{\be1}$, then the function
\lb$t\mapsto \f(U_t)$
has second derivative and
\bay
\label{vpu}
\frac{d^2}{ds^2}\big(\f(U_s)\big)\Big|_{s=0}=-2\left(\iiint(\dg^2\f)(\l,\mu,\nu)
\,dE_U(\l)A\,dE_U(\mu)A\,dE_U(\nu)\right)U^2.
\ey
\end{thm}

Note that by Theorem \ref{tens}, the right-hand side of \rf{vpu} makes sense and determines a bounded
linear operator.

First we prove Theorem \ref{tens} and then we deduce from it
Theorem \ref{uni}.

\medskip

{\bf Proof of Theorem \ref{tens}.} It is easy to see that
\bay
\label{d2}
(\dg^2\f)(z_1,z_2,z_3)=\sum_{i,j,k\ge0}\hat\f(i+j+k+2)z_1^iz_2^jz_3^k+
\sum_{i,j,k\le0}\hat\f(i+j+k-2)z_1^iz_2^jz_3^k,
\ey
where $\hat\f(n)$ is the $n$th Fourier coefficient of $\f$. We prove that
$$
\sum_{i,j,k\ge0}\hat\f(i+j+k+2)z_1^iz_2^jz_3^k\in C(\T)\hat\otimes C(\T)\hat\otimes C(\T).
$$
The fact that 
$$
\sum_{i,j,k\le0}\hat\f(i+j+k-2)z_1^iz_2^jz_3^k\in C(\T)\hat\otimes C(\T)\hat\otimes C(\T)
$$
can be proved in the same way.
Clearly, we can assume that $\hat\f(j)=0$ for $j<0$.

We have
\begin{align*}
\sum_{i,j,k\ge0}\hat\f(i+j+k+2)z_1^iz_2^jz_3^k&=\sum_{i,j,k\ge0}\a_{ijk}\hat\f(i+j+k+2)z_1^iz_2^jz_3^k\\
&+\sum_{i,j,k\ge0}\b_{ijk}\hat\f(i+j+k+2)z_1^iz_2^jz_3^k\\
&+\sum_{i,j,k\ge0}\g_{ijk}\hat\f(i+j+k+2)z_1^iz_2^jz_3^k,
\end{align*}
where
$$
\a_{ijk}=\left\{\begin{array}{ll}\frac13,&i=j=k=0,\\[.2cm]
\frac{i}{i+j+k},& i+j+k\ne0;
\end{array}\right.
$$

\medskip

$$
\b_{ijk}=\left\{\begin{array}{ll}\frac13,&i=j=k=0,\\[.2cm]
\frac{j}{i+j+k},& i+j+k\ne0;
\end{array}\right.
$$
and
$$
\g_{ijk}=\left\{\begin{array}{ll}\frac13,&i=j=k=0,\\[.2cm]
\frac{k}{i+j+k},& i+j+k\ne0.
\end{array}\right.
$$
Clearly, it suffices to show that
\bay
\label{nep}
\sum_{i,j,k\ge0}\a_{ijk}\hat\f(i+j+k+2)z_1^iz_2^jz_3^k\in C(\T)\hat\otimes C(\T)\hat\otimes C(\T).
\ey
It is easy to see that
$$
\sum_{i,j,k\ge0}\a_{ijk}\hat\f(i+j+k+2)z_1^iz_2^jz_3^k=
\sum_{j.k\ge0}\left(\Big(\big((S^*)^{j+k+2}\f\big)*\sum_{i\ge0}\a_{i+j+k}z^i\Big)(z_1)\right)z_2^jz_3^k,
$$
where $S^*$ is backward shift, i.e., $(S^*)^k\f=\pp_+\bar z^k\f$ ($\pp_+$
is the orthogonal projection from $L^2$ onto the Hardy class $H^2$).
Thus
$$
\left\|\sum_{i,j,k\ge0}\a_{ijk}\hat\f(i+j+k+2)z_1^iz_2^jz_3^k\right\|_{L^\be\hat\otimes
L^\be\hat\otimes L^\be}\!\!\!\le
\sum_{j,k\ge0}\left\|\big((S^*)^{j+k+2}\f\big)*\sum_{i\ge0}\a_{i+j+k}z^i\right\|_{L^\be}.
$$
Put
$$
Q_m(z)=\sum_{i\ge m}\frac{i-m}{i}z^i,\quad m>0,\quad \mbox{and}\quad
Q_0(z)=\frac13+\sum_{i\ge1}z^i.
$$
Then it is easy to see that
$$
\left\|\big((S^*)^{j+k+2}\f\big)*\sum_{i\ge0}\a_{i+j+k}z^i\right\|_{L^\be}=\|\psi*Q_{j+k}\|_{L^\be},
$$
where $\psi=(S^*)^2\f$,
and so 
\begin{align*}
\left\|\sum_{i,j,k\ge0}\a_{ijk}\hat\f(i+j+k+2)z_1^iz_2^jz_3^k\right\|_{L^\be\hat\otimes
L^\be\hat\otimes L^\be}&\le\sum_{j,k\ge0}\|\psi*Q_{j+k}\|_{L^\be}\\
&=\sum_{m\ge0}(m+1)\|\psi*Q_m\|_{L^\be}.
\end{align*}

Consider the function $r$ on $\R$ defined by
$$
r(x)=\left\{\begin{array}{ll}1,&|x|\le1,\\[.2cm]\frac1x,&|x|\ge1.
\end{array}\right.
$$
It is easy to see that the Fourier transform $\F r$ of $h$ belongs to $L^1(\R)$.
Define the functions $R_n$, $n\ge1$, on $\T$ by
$$
R_n(\z)=\sum_{k\in\Z}r\left(\frac kn\right)\z^k.
$$

\begin{lem}
\label{Hn}
$$
\|R_n\|_{L^1}\le\const.
$$
\end{lem}

\Pf For $N>0$ consider the function $\xi_N$ defined by
$$
\xi_N(x)=\left\{\begin{array}{ll}1,&|x|\le
N,\\[.2cm]\frac{2N-|x|}N,&N\le|x|\le 2N,\\[.2cm]0,&|x|\ge2N.
\end{array}\right.
$$
It is easy to see that $\F\xi_N\in L^1(\R)$ and
$\|\F\xi_N\|_{L^1(\R)}$ does not depend on $N$.
Let 
$$
R_{N,n}(\z)=\sum_{k\in\Z}r\left(\frac kn\right)\xi_N\left(\frac kn\right)\z^k,\quad\z\in\T.
$$
It was proved in Lemma 2 of \cite{Pe1} that $\|R_{N,n}\|_{L^1}\le\|\F(r\xi_N)\|_{L^1(\R)}$.
Since 
$$
\|\F(r\xi_N)\|_{L^1(\R)}\le\|\F r\|_{L^1(\R)}\|\F\xi_N\|_{L^1(\R)}=\const,
$$
it follows that the $L^1$-norms of $R_{N,n}$ are uniformly bounded. The result
follows from the obvious fact that 
$$
\lim_{N\to\be}\|R_n-R_{N,n}\|_{L^2}=0. \quad \bl
$$

\medskip

Let us complete the proof of Theorem \ref{tens}.

For $f\in L^\be$, we have
$$
\|f*Q_m\|_{L^\be}=\|f-f*R_m\|_{L^\be}\le\|f\|_{L^\be}+\|f*R_m\|_{L^\be}\le\const\|f\|_{L^\be}.
$$

Thus
\begin{align*}
\sum_{m\ge0}(m+1)\|\psi*Q_m\|_{L^\be}&=
\sum_{m\ge0}(m+1)\left\|\sum_{n\ge0}\psi*W_n*Q_m\right\|_{L^\be}\\[.2cm]
&\le\sum_{m,n\ge0}(m+1)\|\psi*W_n*Q_m\|_{L^\be}\\[.1cm]
&=\sum_{n\ge0}\,\,\,\sum_{0\le m\le2^{n+1}}(m+1)\|\psi*W_n*Q_m\|_{L^\be}\\[.1cm]
&\le\const\sum_{n\ge0}\,\,\,\sum_{0\le m\le2^{n+1}}(m+1)\|\psi*W_n\|_{L^\be}\\[.1cm]
&\le\const\sum_{n\ge0}2^{2n}\|\psi*W_n\|_{L^\be}\le\const\|\psi\|_{B_{\be1}^2},
\end{align*}
where the $W_n$ are defined in \S 3.

This proves that 
$$
\sum_{i,j,k\ge0}\a_{ijk}\hat\f(i+j+k+2)z_1^iz_2^jz_3^k\in L^\be\hat\otimes C(\T)\hat\otimes C(\T)
$$
and
\bay
\label{ots}
\left\|\sum_{i,j,k\ge0}\a_{ijk}\hat\f(i+j+k+2)z_1^iz_2^jz_3^k\right\|
_{L^\be\hat\otimes C(\T)\hat\otimes C(\T)}\le\const\|\f\|_{B_{\be1}^2}.
\ey

To prove \rf{nep}, it suffices to represent $\f$ as
$$
\f=\sum_{n\ge0}\f*W_n.
$$
Then we can apply the above reasoning to each polynomial
$\f*W_n$. Since 
$$
\Big(\big((S^*)^{j+k+2}\f*W_n\big)*\sum_{i\ge0}\a_{i+j+k}z^i\Big)
$$ 
is obviously a polynomial, the above reasoning shows that 
$$
\sum_{i,j,k\ge0}\a_{ijk}\widehat{\f*W_n}(i+j+k+2)z_1^iz_2^jz_3^k\in 
C(\T)\hat\otimes C(\T)\hat\otimes C(\T)
$$
and by \rf{ots},
\begin{align*}
\left\|\sum_{i,j,k\ge0}\a_{ijk}\widehat{\f*W_n}(i+j+k+2)z_1^iz_2^jz_3^k\right\|
_{C(\T)\hat\otimes C(\T)\hat\otimes C(\T)}
&\le\const\|\f*W_n\|_{B_{\be1}^2}\\[.3cm]
&\le\const2^{2n}\|\f*W_n\|_{L^\be}.
\end{align*}

The result follows now from the fact that
$$
\sum_{n\ge0}2^{2n}\|\f*W_n\|_{L^\be}\le\const\|\f\|_{B_{\be1}^2}
$$
(see \S 3). $\bl$

Now we are ready to prove Theorem \ref{uni}. 

\medskip

{\bf Proof of Theorem \ref{uni}.} It follows from the definition of the second order divided
difference (see \S 1) that
\bay
\label{div}
(\mu-\nu)(\dg^2\f)(\l,\mu,\nu)=(\dg\f)(\l,\mu)-(\dg\f)(\l,\nu).
\ey

By \rf{vp}, we have
\begin{align*}
&\frac1t\left(\frac{d}{ds}\big(\f(U_s)\big)\Big|_{s=t}-\frac{d}{ds}\big(\f(U_s)\big)\Big|_{s=0}\right)
\\[.2cm]
=&\frac{\rm i}t\left(\iint(\dg\f)(\l,\nu)\,dE_{U_t}(\l)A\,dE_{U_t}(\nu)U_t-
\iint(\dg\f)(\mu,\nu)\,dE_U(\mu)A\,dE_{U_t}(\nu)U\right)\\[.2cm]
=&\frac{\rm i}t\left(\iint(\dg\f)(\l,\nu)\,dE_{U_t}(\l)A\,dE_{U_t}(\nu)-
\iint(\dg\f)(\mu,\nu)\,dE_U(\mu)A\,dE_{U_t}(\nu)\right)U_t\\[.2cm]
&+\frac{\rm i}t\left(\iint(\dg\f)(\mu,\nu)\,dE_U(\mu)A\,dE_{U_t}(\nu)U_t
-\iint(\dg\f)(\mu,\nu)\,dE_U(\mu)A\,dE_{U_t}(\nu)U\right)\\[.2cm]
&+\frac{\rm i}t\left(\iint(\dg\f)(\l,\nu)\,dE_U(\l)A\,dE_{U_t}(\nu)-
\iint(\dg\f)(\l,\mu)\,dE_U(\l)A\,dE_U(\mu)\right)U\\\\[.2cm]
=&\frac{\rm i}t\iiint(\dg^2\f)(\l,\mu,\nu)
\,dE_{U_t}(\l)(e^{{\rm i}tA}-I)U\,dE_U(\mu)A\,dE_{U_t}(\nu)U_t\\[.2cm]
&+\frac{\rm i}t\left(\iint(\dg\f)(\mu,\nu)\,dE_U(\mu)A\,dE_{U_t}(\nu)U_t
-\iint(\dg\f)(\mu,\nu)\,dE_U(\mu)A\,dE_{U_t}(\nu)U\right)\\[.2cm]
&+\frac{\rm i}t\iiint(\dg^2\f)(\l,\mu,\nu)
\,dE_U(\l)A\,dE_U(\mu)(e^{{\rm i}tA}-I)U\,dE_{U_t}(\nu)U_t
\end{align*}
by \rf{div}.

Since $\lim\limits_{t\to0}\|U_t-U\|=0$, to complete the proof it suffices to show that
\begin{align}
\label{per}
\lim_{t\to0}\,\,\,&\frac1t\iiint(\dg^2\f)(\l,\mu,\nu)
\,dE_{U_t}(\l)(e^{{\rm i}tA}-I)U\,dE_U(\mu)A\,dE_{U_t}(\nu)\nonumber\\[.2cm]&=
{\rm i}\iiint(\dg^2\f)(\l,\mu,\nu)\,dE_U(\l)A\,dE_U(\mu)A\,dE_U(\nu)U,
\end{align}
\medskip
\begin{align}
\label{vto}
\lim_{t\to0}\iint(\dg\f)(\mu,\nu)\,dE_U(\mu)A\,dE_{U_t}(\nu)=\iint(\dg\f)(\mu,\nu)\,dE_U(\mu)A\,dE_U(\nu),
\end{align}
and
\begin{align}
\label{tre}
\lim_{t\to0}\,\,\,&\frac1t\iiint(\dg^2\f)(\l,\mu,\nu)
\,dE_U(\l)A\,dE_U(\mu)(e^{{\rm i}tA}-I)U\,dE_{U_t}(\nu)\nonumber\\[.2cm]
&={\rm i}\iiint(\dg^2\f)(\l,\mu,\nu)\,dE_U(\l)A\,dE_U(\mu)A\,dE_U(\nu)U.
\end{align}

Let us prove \rf{per}. Since $\dg^2\f\in C(\T)\hat\otimes C(\T)\hat\otimes C(\T)$, it suffices to show
that for $f,\,g,\,h\in C(\T)$,
\begin{align}
\label{pred}
\lim_{t\to0}\,\,&\frac1t\iiint f(\l)g(\mu)h(\nu)
\,dE_{U_t}(\l)(e^{{\rm i}tA}-I)U\,dE_U(\mu)A\,dE_{U_t}(\nu)\nonumber\\[.2cm]&=
{\rm i}\iiint f(\l)g(\mu)h(\nu)\,dE_U(\l)A\,dE_U(\mu)A\,dE_U(\nu)U.
\end{align}
We have
\begin{align*}
&\frac1t\iiint f(\l)g(\mu)h(\nu)
\,dE_{U_t}(\l)(e^{{\rm i}tA}-I)U\,dE_U(\mu)A\,dE_{U_t}(\nu)\\[.2cm]
=&f(U_t)\left(\frac1t(e^{{\rm i}tA}-I)U\right)g(U)Ah(U_t)
\end{align*}
and
\begin{align*}
\iiint f(\l)g(\mu)h(\nu)\,dE_U(\l)A\,dE_U(\mu)A\,dE_U(\nu)U
=f(U)Ag(U)Ah(U)U.
\end{align*}
Since $f$ and $h$ are in $C(\T)$, it follows that
$$
\lim_{t\to0}\|f(U_t)-f(U)\|=\lim_{t\to0}\|h(U_t)-h(U)\|=0
$$
(it suffices to prove this for trigonometric polynomials $f$ and $h$ which is evident).
This together with the obvious fact
$$
\lim_{t\to0}\left(\frac1t(e^{{\rm i}tA}-I)\right)={\rm i}A
$$
proves \rf{pred} which in turn implies \rf{per}.

The proof of \rf{tre} is similar. To prove \rf{vto}, we observe that $B_{\be1}^2\subset B_{\be1}^1$
and use the fact that $\dg\f\in C(\T)\hat\otimes C(\T)$ (this was proved in \cite{Pe1}). Again, it
suffices to prove that for $f,g \in C(\T)$,
$$
\lim_{t\to0}\iint f(\mu)g(\nu)\,dE_U(\mu)A\,dE_{U_t}(\nu)=\iint f(\mu)g(\nu)\,dE_U(\mu)A\,dE_U(\nu)
$$
which follows from the obvious equality:
$$
\lim_{t\to0}\|g(U_t)-g(U)\|=0.\quad\bl
$$

The proofs of Theorems \ref{tens} and \ref{uni} given above generalize easily to the case of higher
derivatives.

\begin{thm}
\label{mtens}
Let $m$ be a positive integer.
If $\f\in B_{\be1}^m$, then 
$$
\dg^m\f\in\underbrace{C(\T)\hat\otimes\cdots\hat\otimes C(\T)}_{m+1}.
$$
\end{thm}

\begin{thm}
\label{muni}
Let $m$ be a positive integer and 
let $\f$ be a function in the Besov class $B^m_{\be1}$, then the function
$t\mapsto \f(U_t)$
has $m$th derivative and
\begin{align*}
&\frac{d^m}{ds^m}\big(\f(U_s)\big)\Big|_{s=0}\\[.2cm]
=&{\rm i}^mm!\left(\underbrace{\int\cdots\int}_{m+1}(\dg^m\f)(\l_1,\cdots,\l_{m+1})
\,dE_U(\l_1)A\cdots A\,dE_U(\l_{m+1})\right)U^m.
\end{align*}
\end{thm}

\

\section{\bf The case of self-adjoint operators}
\setcounter{equation}{0}

\

In this section we consider the problem of the existence of higher derivatives of the function
$$
t\mapsto \f(A_t)=\f(A+tK)
$$ 
Here $A$ is a self-adjoint operator (not necessarily
bounded), $K$ is a bounded self-adjoint operator, and $A_t\df A+tK$.

In \cite{Pe2} it was shown that if $\f\in B_{\be1}^1(\R)$, then 
$\dg\f\in {\frak B}(\R)\hat\otimes_{\rm i}{\frak B}(\R)$, where 
${\frak B}(\R)$ is the space of bounded Borel functions on $\R$ equipped with the $\sup$-norm,
and so
\bay
\label{ne}
\|\f(A+K)-\f(A)\|\le\const\|\f\|_{B_{\be1}^1}\|K\|.
\ey

In fact, the construction given in \cite{Pe2} shows that for $\f\in B_{\be1}^1(\R)$, the function
$t\mapsto \f(A+tK)$ 
is differentiable and
\bay
\label{fd}
\frac{d}{ds}\big(\f(A_s)\big)\Big|_{s=0}=\iint(\dg\f)(\l,\mu)\,dE_A(\l)K\,dE_A(\mu).
\ey
For completeness, we show briefly how to deduce \rf{fd} from the construction given in \cite{Pe2}. 
We are going to give a detailed proof in the case of higher derivatives.

We need the following notion.

\medskip

{\bf Definition.} A continuous function $\f$ on $\R$ is called {\it operator continuous} if
$$
\lim_{s\to0}\|\f(A+sK)-\f(A)\|=0
$$
for any self-adjoint operator $A$ and any bounded self-adjoint operator $K$. 

\medskip

It follows from \rf{ne} that functions in $B_{\be1}^1(\R)$ are operator continuous.
It is also easy to see that the product of two bounded operator continuous functions is
operator continuous.

\medskip

{\bf Proof of \rf {fd}.} The construction given in
\cite{Pe2} shows that if $\f\in B_{\be1}^1(\R)$, then $\dg\f$ admits a representation
$$
(\dg\f)(\l,\mu)=\int_Q f(\l,x)g(\mu,x)\,d\s(x),
$$
where $(Q,\s)$ is a measure space, $f$ and $g$ are measurable functions on $\R\times Q$ such that
$$
\int_Q\|f_x\|_{{\frak B}(\R)}\|g_x\|_{{\frak B}(\R)}\,d\s(x)<\be,
$$
and for almost all \mbox{$x\in Q$},  and $f_x$ and $g_x$ are operator continuous functions 
where \lb$f_x(\l)\df f(\l,x)$ and $g_x(\mu)\df g(\mu,x)$. Indeed, it is very easy to verify that
the functions $f_x$ and $g_x$ constructed in \cite{Pe2} are products of bounded functions in
$B_{\be1}^1(\R)$.

By \rf{BSF}, we have
\begin{align*}
\frac1s\big(\f(A_s)-\f(A)\big)&=\frac1s\iint(\dg\f)(\l,\mu)\,dE_{A_s}(\l)sK\,dE_A(\mu)\\[.2cm]
&=\int_Q f_x(A_s)Kg_x(A)\,d\s(x).
\end{align*}
Since $f_x$ is operator continuous, we have
$$
\lim_{s\to0}\|f_x(A_s)-f_x(A)\|=0.
$$
It follows that
\begin{align*}
&\left\|\int_Q f_x(A_s)Kg_x(A)\,d\s(x)-\int_Q f_x(A)Kg_x(A)\,d\s(x)\right\|\\[.2cm]
&\le\|K\|\int_Q\|f_x(A_s)-f_x(A)\|\cdot\|g_x(A)\|\,d\s(x)
\to0, \quad\mbox{as}\quad s\to0,
\end{align*}
which implies \rf{fd}. $\bl$

Consider first the problem of the existence of the second operator derivative. First we prove that
if $f\in B_{\be1}^2(\R)$, then 
$\dg^2\f\in{\frak B}(\R)\hat\otimes_{\rm i}{\frak B}(\R)\hat\otimes_{\rm i}{\frak B}(\R)$.
Actually, to prove the existence of the second derivative, we need the following slightly stronger
result.

\begin{thm}
\label{d2f}
Let $\f\in B_{\be1}^2(\R)$. Then there exist a measure space $(Q,\s)$ and measurable functions
$f,\,g$, and $h$ on $\R\times Q$ such that 
\bay
\label{ipr}
(\dg^2\f)(\l,\mu,\nu)=\int_Q f(\l,x)g(\mu,x)h(\nu,x)\,d\s(x),
\ey
$f_x,\,g_x$, and $h_x$ are operator continuous functions for almost all $x\in Q$, and
\bay
\label{tn}
\int_Q\|f_x\|_{{\frak B}(\R)}\|g_x\|_{{\frak B}(\R)}\|h_x\|_{{\frak B}(\R)}\,d\s(x)\le
\const\|f\|_{B_{\be1}^2(\R)}.
\ey
\end{thm}

As before, $f_x(\l)=f(\l,x)$, $g_x(\mu)=g(\mu,x)$, and $h_x(\nu)=h(\nu,x)$.

Theorem \ref{d2f} will be used to prove the main result of this section.

\begin{thm}
\label{vps}
Suppose that $A$ is a self-adjoint operator, $K$ is a bounded self-adjoint operator.
If $\f\in B_{\be1}^2(\R)\bigcap B_{\be1}^1(\R)$, then the function
$s\mapsto\f(A_s)$
has second derivative that is a bounded operator and
\bay
\label{vt}
\frac{d^2}{ds^2}\big(\f(A_s)\big)\Big|_{s=0}=
2\iint(\dg^2\f)(\l,\mu,\nu)\,dE_A(\l)K\,dE_A(\mu)K\,dE_A(\nu).
\ey
\end{thm} 

Note that by Theorem \ref{d2f}, the right-hand side of \rf{vt} makes sense and is a bounded linear
operator.

For $t>0$ and a function $f$, we define $\cS^*_t f$ by
$$
\big(\F(\cS^*_t f)\big)(s)=\left\{\begin{array}{ll}(\F f)(s-t),&t\le s,\\[.2cm]
0,&t>s.
\end{array}\right.
$$
We also define the distributions $q_t$ and $r_t$, $t>0$, by
$$
(\F q_t)(s)=\left\{\begin{array}{ll}\frac{s}{s+t},&s\ge0,\\[.2cm]0,&s<0,
\end{array}\right.
$$
and
$$
(\F r_t)(s)=\left\{\begin{array}{ll}1,&|s|\le t,\\[.2cm]\frac{t}{s},&|s|>t.
\end{array}\right.
$$
It is easy to see that $r_t\in L^1(\R)$ (see \S 4) and $\|r_t\|_{L^1(\R)}$ does not depend on $t$.

To prove Theorem \ref{d2f}, we need the following lemma.

\begin{lem}
\label{2m}
Let $M>0$ and let $\f$ be a bounded function on $\R$ such that \lb$\supp\F\f\subset[M/2,2M]$.
Then 
\begin{align}
\label{tri}
(\dg^2\f)(\l,\mu,\nu)=&
-\iint\limits_{\R_+\times\R_+}\big((\cS_{t+u}^*\f)*q_{t+u}\big)(\l)e^{{\rm i}t\mu}e^{{\rm i}u\nu}
\,dt\,du\nonumber\\[.2cm]
&-\iint\limits_{\R_+\times\R_+}\big((\cS_{s+u}^*\f)*q_{s+u}\big)(\mu)e^{{\rm i}s\l}e^{{\rm i}u\nu}
\,ds\,du\nonumber\\[.2cm]
&-\iint\limits_{\R_+\times\R_+}\big((\cS_{s+t}^*\f)*q_{s+t}\big)(\nu)e^{{\rm i}s\l}e^{{\rm i}t\mu}
\,ds\,dt.
\end{align}
\end{lem}

\Pf Let us first assume that $\F\f\in L^1(\R)$. We have
\begin{align*}
&\iint\limits_{\R_+\times\R_+}\big((\cS_{t+u}^*\f)*q_{t+u}\big)(\l)e^{{\rm i}t\mu}e^{{\rm i}u\nu}
\,d\mu\,d\nu\\[.3cm]
=&\iiint\limits_{\R_+\times\R_+\times\R_+}
(\F\f)(s+t+u)\frac{s}{s+t+u}e^{{\rm i}s\l}e^{{\rm i}t\mu}e^{{\rm i}u\nu}\,ds\,dt\,du.
\end{align*}
We can write similar representations for the other two terms on the right-hand side of \rf{tri}, 
take their sum and
reduce \rf{tri} to the verification of the following identity:
$$
(\dg^2\f)(\l,\mu,\nu)=
-\iiint\limits_{\R_+\times\R_+\times\R_+}
(\F\f)(s+t+u)e^{{\rm i}s\l}e^{{\rm i}t\mu}e^{{\rm i}u\nu}\,ds\,dt\,du.
$$
This identity can be verified elementarily by making the substitution $a=s+t+u$, $b=t+u$, and $c=u$. 

Consider now the general case, i.e., $\f\in L^\be(\R)$ and $\supp\F\f\subset[M/2,2M]$. Consider a
smooth function $\o$ on $\R$ such that $\o\ge0$, $\supp\o\subset[-1,1]$, and $\|\o\|_{L^1(\R)}=1$.
For $\e>0$ we put $\o_\e(x)=\o(x/\e)/\e$ and
define the function $\f_\e$ by $\F\f_\e=(\F\f)*\o_\e$.
Clearly, 
$$
\F\f_\e\in L^1(\R),\quad\lim_{\e\to0}\|\f_\e\|_{L^\be(\R)}=\|\f\|_{L^\be(\R)},
$$
and
$$
\lim_{\e\to0}\f_\e(x)=\f(x)\quad\mbox{for almost all}\quad x\in\R.
$$

Since we have already proved that \rf{tri} holds for $\f_\e$ in place of $\f$, the result follows by
passing to the limit as $\e\to\be$. $\bl$

\medskip

{\bf Proof of Theorem \ref{d2f}.} Suppose that $\supp\F\f\subset[M/2,2M]$.
Let us show that each summand on the right-hand side of
\rf{2m} admits a desired representation. Clearly, it suffices to do it for the first summand.
Put
$$
\psi(\l,\mu,\nu)=
\iint\limits_{\R_+\times\R_+}
\big((\cS_{t+u}^*\f)*q_{t+u}\big)(\l)e^{{\rm i}t\mu}e^{{\rm i}u\nu}\,dt\,du
=\iint\limits_{\R_+\times\R_+}f_{t+u}(\l)g_t(\mu)h_u(\nu)\,dt\,du,
$$
where
$$
f_v(\l)=\big((\cS_v^*\f)*q_v\big)(\l),\quad g_t(\mu)=e^{{\rm i}t\mu},\quad
\mbox{and}\quad h_u(\nu)=e^{{\rm i}u\nu}.
$$
Clearly, $\|g_t\|_{{\frak B}(\R)}=1$ and $\|h_u\|_{{\frak B}(\R)}=1$. Since
$$
\|f_v\|_{{\frak B}(\R)}=\|f_v\|_{L^\be}=\|\f-\f*r_v\|_{L^\be}
\le
\left\{\begin{array}{ll}(1+\|r_v\|_{L^1})\|\f\|_{L^\be},&v\le2M,\\[.2cm]0,&v>2M,
\end{array}\right.
$$
we have
$$
\|\psi\|_{{\frak B}(\R)\hat\otimes_{\rm i}{\frak B}(\R)\hat\otimes_{\rm i}{\frak B}(\R)}\le
\const\|\f\|_{L^\be}\iint\limits_{t,u>0,t+u\le 2M}\,dt\,du
\le\cs M^2\|\f\|_{L^\be}.
$$

In the same way we can treat the case when $\supp\F\f\subset[-2M,-M/2]$. If $\f$ is a polynomial
of degree at most 2, the result is trivial. 

Let now $\f\in B_{\be1}^1(\R)$ and
$$
\f=\sum_{n\in\Z}\f*W_n+\sum_{n\in\Z}\f*W^\#.
$$
It follows from the above estimate that
$$
\|\dg^2\f\|_{{\frak B}(\R)\hat\otimes_{\rm i}{\frak B}(\R)\hat\otimes_{\rm i}{\frak B}(\R)}\le\
\const\left(\sum_{n\in\Z}2^{2n}\|\f*W_n\|_{L^\be}+
\sum_{n\in\Z}2^{2n}\|\f*W^\#_n\|_{L^\be}\right).
$$

To complete the proof of Theorem \ref{d2f}, we observe that the functions $\l\mapsto e^{{\rm i}tv}$
are operator continuous, because they belong to $B^1_{\be1}(\R)$. On the other hand, it is easy to see
that if $\supp\f\subset[M/2.2M]$, then
the function $(\cS_v^*\f)*q_v$ is the product of $e^{{\rm i}tv}$ and a bounded function
in $B^1_{\be1}(\R)$. $\bl$

\medskip

To prove \rf{vps}, we need the following lemma.

\begin{lem}
\label{ol}
Let $A$ be a self-adjoint operator and let $K$ be a bounded self-adjoint operator.
Suppose that $\f$ is a function on $\R$ such that $\dg\f\in L^\be(\R)\pt L^\be(\R)$
and $\dg^2\f\in L^\be(\R)\pt L^\be(\R)\pt L^\be(\R)$. Then
\begin{align*}
&\iint(\dg\f)(\l,\mu)\,dE_{A+K}(\l)K\,dE_{A+K}(\mu)-
\iint(\dg\f)(\l,\nu)\,dE_{A+K}(\l)K\,dE_A(\nu)\\[.2cm]
=&\iiint(\dg^2\f)(\l,\mu,\nu)\,dE_{A+K}(\l)K\,dE_{A+K}(\mu)K\,dE_A(\nu).
\end{align*}
\end{lem}

\Pf Put
$$
P_n=E_A\big([-n,n]\big),\quad Q_n=E_{A+K}\big([-n,n]\big),\quad
A_{[n]}=P_nA,\quad\mbox{and}\quad B_{[n]}=Q_n(A+K).
$$
We have
\begin{align*}
&\iint(\dg\f)(\l,\mu)\,dE_{A+K}(\l)K\,dE_{A+K}(\mu)-
\iint(\dg\f)(\l,\nu)\,dE_{A+K}(\l)K\,dE_A(\nu)\\[.2cm]
=&\iiint(\dg\f)(\l,\mu)\,dE_{A+K}(\l)K\,dE_{A+K}(\mu)\,dE_A(\nu)\\[.2cm]
&-\iiint(\dg\f)(\l,\nu)\,dE_{A+K}(\l)K\,dE_{A+K}(\mu)\,dE_A(\nu).
\end{align*}

Thus
\begin{align*}
&Q_n\left(\iint(\dg\f)(\l,\mu)\,dE_{A+K}(\l)K\,dE_{A+K}(\mu)-
\iint(\dg\f)(\l,\mu)\,dE_{A+K}(\l)K\,dE_A(\mu)\right)P_n\\[.2cm]
=&\int\limits_{-n}^n\int\limits_{-n}^n\int\limits_{-n}^n
(\dg\f)(\l,\mu)\,dE_{A+K}(\l)K\,dE_{A+K}(\mu)\,dE_A(\nu)\\[.2cm]
&-\int\limits_{-n}^n\int\limits_{-n}^n\int\limits_{-n}^n
(\dg\f)(\l,\nu)\,dE_{A+K}(\l)K\,dE_{A+K}(\mu)\,dE_A(\nu)\\[.2cm]
=&\iiint
(\mu-\nu)(\dg^2\f)(\l,\mu,\nu)\,dE_{B_{[n]}}(\l)K\,dE_{B_{[n]}}(\mu)\,dE_{A_{[n]}}(\nu),
\end{align*}
since
$$
(\dg\f)(\l,\mu)-(\dg\f)(\l,\nu)=(\mu-\nu)(\dg^2\f)(\l,\mu,\nu).
$$

On the other hand,
\begin{align*}
&Q_n\left(\iiint(\dg^2\f)(\l,\mu,\nu)\,dE_{A+K}(\l)K\,dE_{A+K}(\mu)K\,dE_A(\nu)\right)P_n\\[.2cm]
=&\int\limits_{-n}^n\int\limits_{-n}^n\int\limits_{-n}^n
(\dg^2\f)(\l,\mu,\nu)\,dE_{A+K}(\l)K\,dE_{A+K}(\mu)\Big((A+K)-A\Big)\,dE_A(\nu)\\[.2cm]
=&\int\limits_{-n}^n\int\limits_{-n}^n\int\limits_{-n}^n
(\dg^2\f)(\l,\mu,\nu)\,dE_{A+K}(\l)K\,dE_{A+K}(\mu)Q_n\Big((A+K)-A\Big)P_n\,dE_A(\nu)\\[.2cm]
=&\int\limits_{-n}^n\int\limits_{-n}^n\int\limits_{-n}^n
(\dg^2\f)(\l,\mu,\nu)\,dE_{A+K}(\l)K\,dE_{A+K}(\mu)(B_{[n]}-A_{[n]})\,dE_A(\nu)\\[.2cm]
=&\iiint(\dg^2\f)(\l,\mu,\nu)
\,dE_{B_{[n]}}(\l)K\,dE_{B_{[n]}}(\mu)(B_{[n]}-A_{[n]})\,dE_{A_{[n]}}(\nu)\\[.2cm]
\end{align*}
It is easy to see that this is equal to
\begin{align*}
&\iiint(\dg^2\f)(\l,\mu,\nu)\,dE_{B_{[n]}}(\l)K\,dE_{B_{[n]}}(\mu)B_{[n]}\,dE_{A_{[n]}}(\nu)\\[.2cm]
&-\iiint(\dg^2\f)(\l,\mu,\nu)\,dE_{B_{[n]}}(\l)K\,dE_{B_{[n]}}(\mu)A_{[n]}\,dE_{A_{[n]}}(\nu)\\[.2cm]
=&\iiint\mu(\dg^2\f)(\l,\mu,\nu)\,dE_{B_{[n]}}(\l)K\,dE_{B_{[n]}}(\mu)\,dE_{A_{[n]}}(\nu)\\[.2cm]
&-\iiint\nu(\dg^2\f)(\l,\mu,\nu)\,dE_{B_{[n]}}(\l)K\,dE_{B_{[n]}}(\mu)\,dE_{A_{[n]}}(\nu)\\[.2cm]
=&\iiint(\mu-\nu)(\dg^2\f)(\l,\mu,\nu)\,dE_{B_{[n]}}(\l)K\,dE_{B_{[n]}}(\mu)\,dE_{A_{[n]}}(\nu).
\end{align*}

The result follows now from the fact that 
$$
\lim_{n\to\be}P_n=\lim_{n\to\be}Q_n=I
$$
in the strong operator topology. $\bl$

\medskip

{\bf Proof of Theorem \ref{vps}.} It follows from Lemma \ref{ol} that
\begin{align*}
&\frac1t\left(\iint(\dg\f)(\l,\mu)\,dE_{A_t}(\l)K\,dE_{A_t}(\mu)-
\iint(\dg\f)(\l,\nu)\,dE_{A_t}(\l)K\,dE_A(\nu)\right)\\
=&\iiint(\dg^2\f)(\l,\mu,\nu)\,dE_{A_t}(\l)K\,dE_{A_t}(\mu)K\,dE_A(\nu).
\end{align*}

%

Similarly,
\begin{align*}
&\frac1t\left(\iint(\dg\f)(\l,\nu)\,dE_{A_t}(\l)K\,dE_A(\nu)
-\iint(\dg\f)(\mu,\nu)\,dE_A(\mu)K\,dE_A(\nu)\right)\\
=&\iiint(\dg^2\f)(\l,\mu,\nu)\,dE_{A_t}(\l)K\,dE_A(\mu)K\,dE_{A_t}(\nu).
\end{align*}

Thus
\begin{align*}
\frac1t\left(\frac{d}{ds}\f(A_s)\Big|_{s=t}-\frac{d}{ds}\f(A_s)\Big|_{s=0}\right)
&=\iiint(\dg^2\f)(\l,\mu,\nu)\,dE_{A_t}(\l)K\,dE_{A_t}(\mu)K\,dE_A(\nu)\\
&+\iiint(\dg^2\f)(\l,\mu,\nu)\,dE_{A_t}(\l)K\,dE_A(\mu)K\,dE_A(\nu).
\end{align*}

The fact that
\begin{align*}
&\lim_{t\to0}\,\iiint(\dg^2\f)(\l,\mu,\nu)\,dE_{A_t}(\l)K\,dE_{A_t}(\mu)K\,dE_A(\nu)\\[.2cm]
=&
\iiint(\dg^2\f)(\l,\mu,\nu)\,dE_A(\l)K\,dE_A(\mu)K\,dE_A(\nu)
\end{align*}
follows immediately from \rf{ipr} and \rf{tn} and from the fact that the functions
$f_x$, $g_x$, and $h_x$ in \rf{ipr} are operator continuous.

Similarly,
\begin{align*}
&\lim_{t\to0}\iiint(\dg^2\f)(\l,\mu,\nu)\,dE_{A_t}(\l)K\,dE_A(\mu)K\,dE_A(\nu)\\[.2cm]
=&
\iiint(\dg^2\f)(\l,\mu,\nu)\,dE_A(\l)K\,dE_A(\mu)K\,dE_A(\nu),
\end{align*}
which completes the proof. $\bl$

\medskip

{\bf Remark.} In the case of functions on the real line the Besov space 
$B_{\be1}^2(\R)$ is not contained in the space $B_{\be1}^1(\R)$. In the statement of Theorem
\ref{vps} we impose the assumption that $\f\in B_{\be1}^1(\R)$ to ensure that 
the function $t\mapsto\f(A_t)$ has the first derivative. However, we can define the second derivative 
of this function in a slightly different way.

Suppose that $\f\in B_{\be1}^2(\R)$ and
$$
\f=\sum_{n\in\Z}\f*W_n+\sum_{n\in\Z}\f*W_n^\#.
$$
Then the functions $\f_n\df\f*W_n$ and $\f_n^\#\df\f*W_n^\#$ belong to
$B_{\be1}^2(\R)\bigcap B_{\be1}^2(\R)$ and by Theorems \ref{d2f} and \ref{vps}, the series
$$
\sum_{n\in\Z}\frac{d^2}{ds^2}\big(\f_n(A_s)\big)\Big|_{s=0}+
\sum_{n\in\Z}\frac{d^2}{ds^2}\big(\f_n^\#(A_s)\big)\Big|_{s=0}
$$
converges absolutely and we can define the second derivative of the function
$t\mapsto\f(A_t)$ by
$$
\frac{d^2}{ds^2}\big(\f(A_s)\big)\Big|_{s=0}
\df\sum_{n\in\Z}\frac{d^2}{ds^2}\big(\f_n(A_s)\big)\Big|_{s=0}+
\sum_{n\in\Z}\frac{d^2}{ds^2}\big(\f_n^\#(A_s)\big)\Big|_{s=0}.
$$
With this definition the function {\it the function $t\mapsto\f(A_t)$
can possess the second derivative without having the first derivative!}

If $\f(\l)=\l^2$, we can write formally
$$
\frac{d^2}{ds^2}\big(\f(A_s)\big)\Big|_{s=0}=
\frac{d^2}{ds^2}\big(A^2+s(KA+AK)+s^2K^2\big)\Big|_{s=0}=2K^2.
$$
Then formula \rf{vt} holds for an arbitrary $\f\in B_{\be1}^2(\R)$.

\medskip

The proofs of Theorems \ref{d2f} and \ref{vps} given above easily generalize to the case of
derivatives of an arbitrary order.

\begin{thm}
\label{ddf}
Let $m$ be a positive integer and let
$\f\in B_{\be1}^m(\R)$. Then there exist a measure space $(Q,\s)$ and measurable functions
$f_1,\cdots,f_{m+1}$ on $\R\times Q$ such that 
$$
(\dg^m\f)(\l_1,\cdots,\l_{m+1})=\int_Q f_1(\l_1,x)f_2(\l_2,x)\cdots f_{m+1}(\l_{m+1},x)\,d\s(x),
$$
the functions $f_1(\cdot,x),\cdots,f_{m+1}(\cdot,x)$ are operator continuous for almost
all $x\in Q$, and
$$
\int_Q\|f_1(\cdot,x)\|_{{\frak B}(\R)}\cdots\|f_{m+1}(\cdot,x)\|_{{\frak B}(\R)}\,d\s(x)\le
\const\|f\|_{B_{\be1}^m(\R)}.
$$
\end{thm}

\begin{thm}
\label{dps}
Let $m$ be a positive integer.
Suppose that $A$ is a self-adjoint operator, $K$ is a bounded self-adjoint operator.
If $\f\in B_{\be1}^m(\R)\bigcap B_{\be1}^1(\R)$, then the function
\lb$s\mapsto\f(A_s)$
has $m$th derivative that is a bounded operator and
$$
\frac{d^m}{ds^m}\big(\f(A_s)\big)\Big|_{s=0}=
m!\underbrace{\int\cdots\int}_{m+1}(\dg^{m}\f)(\l_1,\cdots,\l_{m+1})
\,dE_A(\l_1)K\cdots K\,dE_A(\l_{m+1}).
$$
\end{thm} 

As in the case $m=2$, we can slightly change the definition of the $m$th derivative so that
for functions $\f\in B_{\be1}^m(\R)$ the function $s\mapsto\f(A_s)$ has $m$th derivative, but does not
have to possess derivatives of orders less than $m$ (see the Remark following the proof of
Theorem \ref{vps}).

\medskip

{\bf Remark.} It is easy to see that in case $A$ is a bounded self-adjoint operator for the existence
of the $m$th derivative of the function $s\mapsto\f(A_s)$, it suffices to assume that $\f$
belongs to $B^m_{\be1}$ locally, i.e., for each finite
interval $I$ there exists a function $\psi$ in $B^m_{\be1}(\R)$ such that $\f\big|I=\psi\big|I$.

\

\

\noindent
\begin{tabular}{p{8cm}p{14cm}}
Department of Mathematics \\
Michigan State University  \\
East Lansing, Michigan 48824\\
USA
\end{tabular}

\end{document}